\documentclass[11pt]{article}
\usepackage[T1]{fontenc}
\usepackage{lmodern}
\usepackage[margin=1in]{geometry}
\usepackage{amsmath,amssymb,amsthm,mathtools}
\usepackage{graphicx,enumitem,microtype,needspace}
\usepackage[font=small,labelfont=bf]{caption}
\usepackage{xcolor}
\usepackage[colorlinks=true,linkcolor=blue!45!black,citecolor=blue!45!black,urlcolor=blue!45!black]{hyperref}
\graphicspath{{starlike_figures/}}
\hypersetup{pdftitle={A starlike 3-dimensional counterexample to the hot spots conjecture},pdfauthor={Rami Atar and Krzysztof Burdzy}}
\setlist[itemize]{leftmargin=1.5em,itemsep=2pt,topsep=3pt}
\setlist[enumerate]{leftmargin=1.7em,itemsep=3pt,topsep=3pt}
\theoremstyle{plain}
\newtheorem{theorem}{Theorem}
\newtheorem{lemma}{Lemma}
\newcommand{\R}{\mathbb R}
\newcommand{\Sph}{S^2}
\newcommand{\B}{B}

\newcommand{\q}{\mathcal E}
\newcommand{\normm}[1]{\|#1\|_m}
\newcommand{\dd}{\,d}
\title{A starlike 3-dimensional counterexample\\to the hot spots conjecture\thanks{This is a first draft. An extended version is expected soon.}}
\author{Rami Atar\thanks{Viterbi Faculty of Electrical and Computer Engineering, Technion--Israel Institute of Technology, Haifa, Israel. E-mail: \texttt{rami@technion.ac.il}.}
\and Krzysztof Burdzy\thanks{Department of Mathematics, University of Washington, Seattle, Washington, USA. E-mail: \texttt{burdzy@uw.edu}.}}
\date{October 14, 2026}
\begin{document}
\maketitle

\begin{abstract}
We construct a bounded Lipschitz domain in $\R^3$ that is starlike
about the origin and whose first positive Neumann eigenvalue is simple. A corresponding eigenfunction,
after a choice of sign, has its unique global maximum at the origin.
The construction uses domains with increasingly fine angular structure
and convergence to a weighted eigenvalue problem whose first positive
eigenfunction is radial.
\end{abstract}

\section{Introduction}\label{intro}

Our main result is the following.

\begin{theorem}\label{main}
There exists a bounded Lipschitz domain $\Omega\subset\R^3$, starlike about $0$, whose first positive Neumann eigenvalue is simple. A corresponding real eigenfunction $u$ has a continuous representative on $\overline\Omega$ and, after a choice of sign, satisfies
\begin{equation}\label{exact}
u(0)>u(x)\qquad\bigl(x\in\overline\Omega\setminus\{0\}\bigr).
\end{equation}
In particular,
\begin{equation}\label{target}
u(0)>\max_{x\in\partial\Omega}u(x).
\end{equation}
\end{theorem}

Rauch's hot spots conjecture \cite{Rauch} asserts that a first nonconstant
Neumann eigenfunction attains its extrema only on the boundary.
Burdzy and Werner \cite{BW} disproved this for a planar domain with two
holes. Bass and Burdzy \cite{BassBurdzy} constructed examples with both
extrema strictly in the interior, and Burdzy \cite{BurdzyOneHole}
obtained such an example with only one hole. Further numerical examples
were given by Kleefeld \cite{Kleefeld}. For an overview, see Burdzy
\cite{BurdzyOverview}.

Positive results for planar domains include those of Ba\~nuelos and
Burdzy \cite{BanuelosBurdzy}, the symmetry results of Jerison and
Nadirashvili \cite{JerisonNadirashvili} and Pascu \cite{Pascu}, the
lip-domain theorem of Atar and Burdzy \cite{AtarBurdzy}, and the triangle
theorem of Judge and Mondal \cite{JudgeMondal}, together with its erratum
\cite{JudgeMondalErratum}. More recently, OpenAI \cite{OpenAI} presented
a proof for smooth simply connected planar domains. Taking that result
as an input, Burdzy \cite{BurdzySmoothness} derived the non-strict
hot-spots property for nonsmooth boundary.

In higher dimensions, Kennedy and Rohleder \cite{KennedyRohleder} proved
positive results for classes extending planar lip domains and symmetric
domains. In contrast, de Dios Pont \cite{deDiosPont} constructed convex
counterexamples in sufficiently high dimensions. De Dios Pont, Hsu and
Taylor \cite{DHT} determined the sharp maximal ratio between the interior
and boundary maxima; their Remark~1.11 also discusses counterexamples
homeomorphic to a ball in dimensions at least three. The present result
shows that failure is possible in dimension three even under the
stronger geometric requirement of starlikeness.

\subsection{Notation}

We write $|x|$ for the Euclidean norm of $x\in\R^3$,
\[
\B(x,r)=\{y\in\R^3:|y-x|<r\},\qquad
\B_r=\B(0,r),\qquad \Sph=\{x\in\R^3:|x|=1\}.
\]
The measure $dS$ is ordinary surface measure on $\Sph$, normalized by $\int_{\Sph}dS=4\pi$. For $U\subset\Sph$, $|U|$ denotes its surface area. The symbols $\overline D$ and $\partial D$ denote the closure and boundary of a set $D$.

A domain $D$ is \emph{starlike about $0$} if $x\in D$ implies $tx\in D$ for every $0\leq t\leq1$. It is \emph{centrally symmetric} if $D=-D$. A Lipschitz domain is one whose boundary is locally the graph of a Lipschitz function in suitable orthogonal coordinates, with the domain on one side of the graph. Corners are allowed.

For an open set $D$, $H^1(D)$ is the space of square-integrable functions with square-integrable first weak derivatives. A Neumann eigenfunction with eigenvalue $\lambda$ is a nonzero $u\in H^1(D)$ satisfying
\[
\int_D\nabla u\cdot\nabla w\dd x
=\lambda\int_D uw\dd x\qquad(w\in H^1(D)).
\]
On a bounded connected Lipschitz domain we number the eigenvalues, with multiplicities, as
\[
0=\lambda_0(D)<\lambda_1(D)\leq\lambda_2(D)\leq\cdots.
\]
Thus $\lambda_1(D)$ is the first nonconstant eigenvalue, also called the second Neumann eigenvalue. All functions and variational spaces in the proof are real unless Fourier series are being used.

For a function on $\B_R$ we write $v(r,\omega)=v(r\omega)$ in polar coordinates. The tangential gradient $\nabla_{\Sph}v(r,\omega)$ always refers to differentiation in $\omega$ on the unit sphere. Matrix transposition is denoted by $\mathsf T$. The constants denoted by $C$ may change from one occurrence to the next; their parameter dependence is specified where needed. The construction uses a large radius $R$, fixed before the small spacing parameter $\delta$ tends to zero.

\subsection{Outline and organization}

The domain contains the unit ball and extends to radius $R$ along an explicitly chosen set of directions. As the angular spacing $\delta$ tends to zero, its Neumann problem converges to a weighted variational problem on $\B_R$. The first positive eigenfunction of that limiting problem is radial and strictly decreasing. Uniform convergence of the actual eigenfunctions transfers the strict difference between the value at the origin and the boundary values. Central symmetry and interior derivative estimates locate the unique maximum exactly at the origin.

Section~\ref{construction} defines the domains. Section~\ref{effective} describes the limiting variational problem and its angular coefficient. Section~\ref{auxiliary} states the six auxiliary lemmas. Using these lemmas, Section~\ref{proof} proves Theorem~\ref{main}. All six lemma proofs are collected in Section~\ref{lemmproofs}.

\section{Construction of the domains}\label{construction}

\subsection{Parameters}

For $R>1$, the spherical annulus between radii $1$ and $R$ is
\[
\mathcal A_{1,R}=\{x\in\R^3:1<|x|<R\}.
\]
The domain will contain only selected directions in this annulus.

Choose $R>1$ and set
\begin{equation}\label{theta}
\alpha=R^{-5/2},\qquad \theta=\frac{\alpha}{4\pi},
\qquad h=1-\sqrt{1-\theta}.
\end{equation}
The number $\theta$ is the limiting fraction of the sphere that will be retained; $\alpha=4\pi\theta$ is its limiting area. The identity $2h-h^2=\theta$ explains the choice of $h$.

The spacing parameter is $\delta>0$. The order of choices is: choose $R$ large enough, \emph{fix it}, and then let $\delta\downarrow0$. Constants in extension and convergence estimates may depend on this fixed $R$ and its associated $\theta$.

\subsection{A subset of the plane: horizontal and vertical strips}\label{planar}

For $t\in\R$, define
\[
\operatorname{dist}(t,\delta\mathbb Z)=\min_{n\in\mathbb Z}|t-\delta n|.
\]
The planar set to be retained is
\begin{equation}\label{planarset}
P_{\delta,h}=\Big\{y=(y_1,y_2)\in\R^2:
\operatorname{dist}(y_1,\delta\mathbb Z)<h\delta/2
\ \hbox{or}\ 
\operatorname{dist}(y_2,\delta\mathbb Z)<h\delta/2\Big\}.
\end{equation}
Thus it is the union of vertical strips centered at $y_1=\delta n$ and horizontal strips centered at $y_2=\delta n$, each of total width $h\delta$.

It is \emph{periodic} in the precise sense that
\[
y\in P_{\delta,h}\quad\Longleftrightarrow\quad
y+\delta a\in P_{\delta,h}\qquad(a\in\mathbb Z^2).
\]
A \emph{period square} is any square $\delta a+[0,\delta]^2$ on which this repeated pattern is represented once.

For $a=(a_1,a_2)\in\mathbb Z^2$, define the closed square
\begin{equation}\label{holes}
H_a=\prod_{\nu=1}^2
\big[a_\nu\delta+h\delta/2,\ (a_\nu+1)\delta-h\delta/2\big].
\end{equation}
Then
\[
\R^2\setminus P_{\delta,h}=\bigcup_{a\in\mathbb Z^2}H_a.
\]
Each $H_a$ has side length $(1-h)\delta$. The retained area in a period square is
\[
\delta^2-(1-h)^2\delta^2=(2h-h^2)\delta^2=\theta\delta^2.
\]

\begin{figure}[!htbp]
\centering\includegraphics[width=\textwidth]{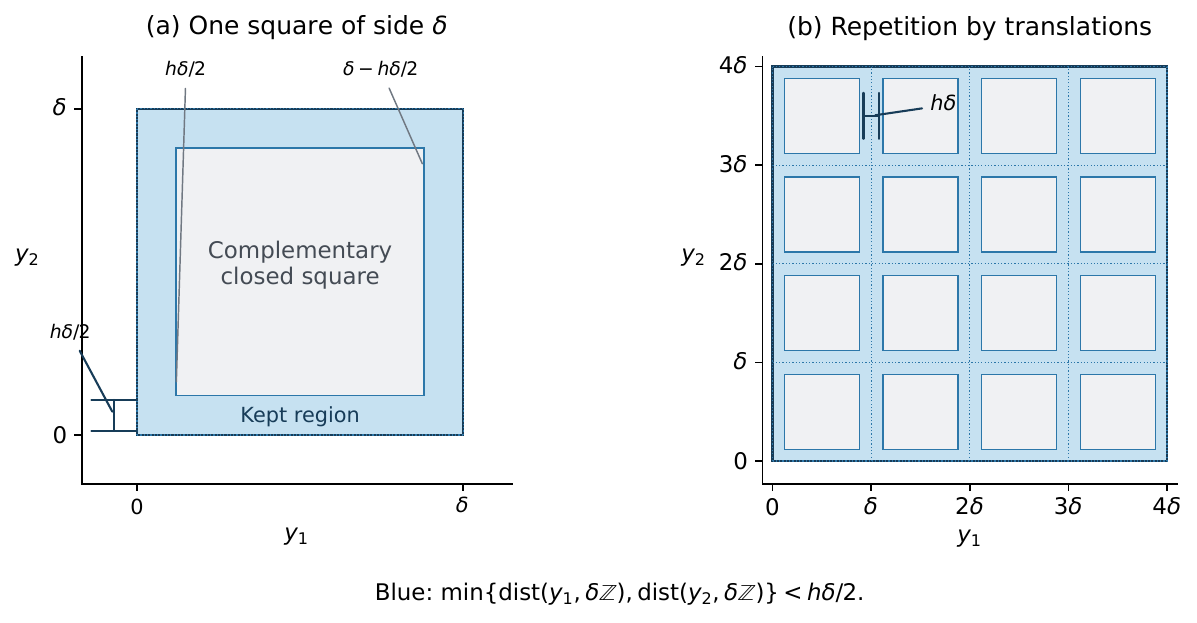}
\caption{The set $P_{\delta,h}$ in one period and in sixteen periods. Blue is retained; grey is excluded, including its boundary. Across adjacent periods the retained strips have total width $h\delta$. For visibility the drawing uses a larger $h$ than the parameter chosen in \eqref{theta}.}
\label{figplanar}
\end{figure}

\subsection{Six explicit coordinate maps for the sphere}\label{coordinates}

Let $e_1,e_2,e_3$ be the standard coordinate vectors. For $i\in\{1,2,3\}$ and $\sigma\in\{-1,1\}$, define the closed spherical region
\begin{equation}\label{faces}
K_{i,\sigma}=\{\omega\in\Sph:
\sigma\omega_i\geq|\omega_j|\ \hbox{for every }j\ne i\}.
\end{equation}
These six regions cover the sphere; their interiors do not overlap. They are the radial projections of the six faces of a cube.

If $j<k$ are the other two indices, define
\begin{equation}\label{psimap}
\Psi_{i,\sigma}(y_1,y_2)
=\frac{\sigma e_i+y_1e_j+y_2e_k}
{\sqrt{1+y_1^2+y_2^2}},
\qquad (y_1,y_2)\in[-1,1]^2.
\end{equation}
It maps the closed coordinate square onto $K_{i,\sigma}$. On the interiors it is smooth, one-to-one, and has the smooth inverse
\[
y_1=\frac{\omega_j}{|\omega_i|},\qquad
y_2=\frac{\omega_k}{|\omega_i|}.
\]

For example,
\[
\Psi_{1,+}(y_1,y_2)=\frac{(1,y_1,y_2)}{\sqrt{1+y_1^2+y_2^2}}.
\]
Define the set of common face boundaries by
\begin{equation}\label{sigmadef}
\Sigma=\bigcup_{i,\sigma}\partial_{\Sph}K_{i,\sigma}.
\end{equation}
Equivalently, $\Sigma$ consists of points where at least two of $|\omega_1|,|\omega_2|,|\omega_3|$ tie for the largest value. It is a finite union of great-circle arcs.

\begin{figure}[!htbp]
\centering\includegraphics[width=\textwidth]{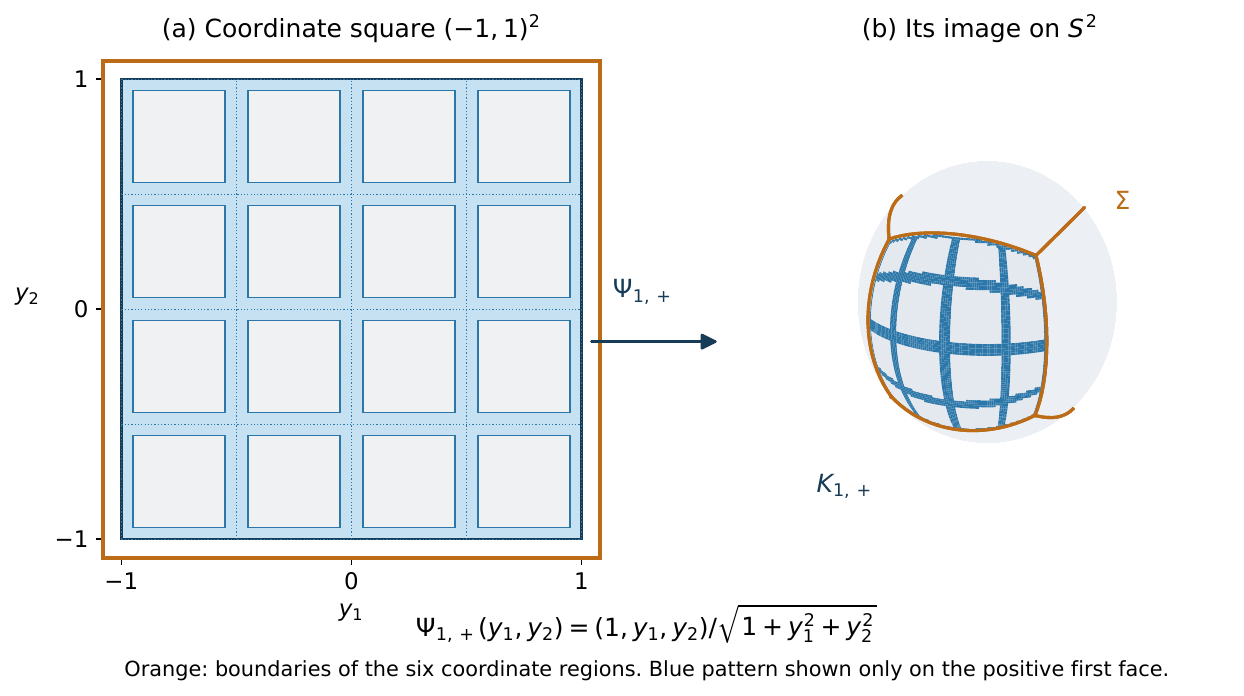}
\caption{The explicit map $\Psi_{1,+}$ sends the square to $K_{1,+}$. Orange curves mark $\Sigma$. Blue shows the image of the planar strips on this one face; the other faces have been left plain to show the geometry. This picture illustrates the coordinate map before the boundary-neighborhood rule in the next subsection is applied.}
\label{figcoords}
\end{figure}

\subsection{Exactly which squares are removed from the sphere}\label{selection}

For $\omega,\eta\in\Sph$, let
\[
d_{\Sph}(\omega,\eta)=\arccos(\omega\cdot\eta),
\qquad d_{\Sph}(\omega,\Sigma)=\min_{\eta\in\Sigma}d_{\Sph}(\omega,\eta).
\]
These are distance along the sphere and distance to its face boundaries. Fix $M=4$ and define the closed neighborhood
\begin{equation}\label{boundaryregion}
N_\delta=\{\omega\in\Sph:d_{\Sph}(\omega,\Sigma)\leq M\delta\}.
\end{equation}
Every point of $N_\delta$ will be retained.

For each face, let $I_{i,\sigma,\delta}$ consist of the indices $a\in\mathbb Z^2$ for which
\begin{equation}\label{indexrule}
H_a\subset(-1,1)^2,
\qquad
\min_{y\in H_a}d_{\Sph}(\Psi_{i,\sigma}(y),\Sigma)>M\delta.
\end{equation}
Only finitely many indices satisfy this rule. Define
\begin{equation}\label{Udef}
U_\delta=\Sph\setminus
\bigcup_{i,\sigma}\ \bigcup_{a\in I_{i,\sigma,\delta}}
\Psi_{i,\sigma}(H_a).
\end{equation}
The set $U_\delta$ consists of the retained directions.

The removal rule can be read as follows. First consider all complementary squares $H_a$ in the coordinate plane. Remove the image of such a square only when the entire square lies inside the coordinate face and its image stays farther than $M\delta$ from $\Sigma$. A square meeting the retained boundary neighborhood is kept in its entirety. Thus no square is cut into a small fragment when two coordinate regions meet.

The removed sets are disjoint closed images of squares. In particular, $U_\delta$ is open. Its connectedness and the remaining geometric estimates are included in Lemma~\ref{geom}.

The construction is symmetric under $\omega\mapsto-\omega$. Indeed,
\[
\Psi_{i,-}(-y)=-\Psi_{i,+}(y),\qquad
-H_a=H_{-a-(1,1)},
\]
and the removal rule is unchanged by this reflection. Hence $U_\delta=-U_\delta$.

\subsection{From directions on the sphere to a domain in three dimensions}\label{radialdomain}

The final definition is
\begin{equation}\label{domain}
\Omega_\delta=\B_1\ \cup\ \{r\omega:0<r<R,\ \omega\in U_\delta\}.
\end{equation}
Equivalently, a nonzero $x$ belongs to $\Omega_\delta$ exactly when
\[
|x|<R\quad\hbox{and}\quad
\big(|x|<1\ \hbox{or}\ x/|x|\in U_\delta\big).
\]
The origin belongs to the entire unit ball. The part outside the unit ball is therefore precisely
\[
\{r\omega:1\leq r<R,\ \omega\in U_\delta\}.
\]

For $1<r<R$,
\[
\Omega_\delta\cap r\Sph=rU_\delta,
\qquad \operatorname{area}(rU_\delta)=r^2|U_\delta|.
\]
For $r<1$ every direction is included. At $r=1$, the openings $\omega\in U_\delta$ are included as well.

\begin{figure}[!htbp]
\centering\includegraphics[width=\textwidth]{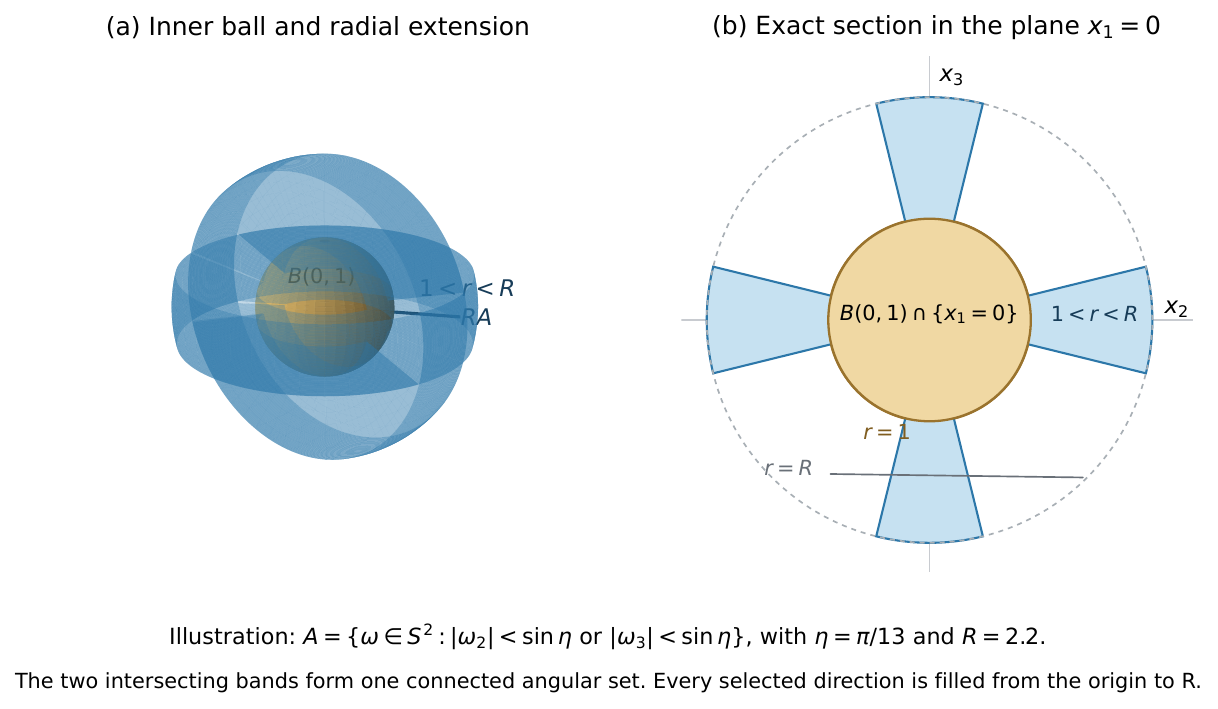}
\caption{The radial definition illustrated with the simpler connected set
$A=\{\omega:|\omega_2|<\sin(\pi/13)\ \hbox{or}\ |\omega_3|<\sin(\pi/13)\}$ and $R=2.2$.
Left: transparent surfaces show $\B_1\cup\{r\omega:0<r<R,\omega\in A\}$.
Right: its exact section in the plane $x_1=0$. This picture explains how angular directions are extended; the set used in the proof is the explicitly defined $U_\delta$ in \eqref{Udef}.}
\label{figdomain}
\end{figure}

\section{The limiting eigenvalue problem}\label{effective}

The actual eigenfunction solves a Neumann problem on $\Omega_\delta$. The limiting problem is a different variational problem on the fixed ball $\B_R$. The precise convergence statement connecting them is Lemma~\ref{hom}.

\subsection{Weighted norm and gradient energy}

Define $m(x)=1$ for $|x|<1$ and $m(x)=\theta$ for $1<|x|<R$. Its values on $|x|=1$ do not affect any integral. For $v\in H^1(\B_R)$,
\begin{equation}\label{mass}
\normm{v}^2=\int_{\B_R}m(x)v(x)^2\dd x
=\int_{\B_1}v^2\dd x
 +\theta\int_1^R\int_{\Sph}v(r,\omega)^2r^2\dd S\dd r.
\end{equation}
The inner product $\langle v,w\rangle_m$ is the same formula with $vw$ instead of $v^2$. The word ``mass'' here denotes this weight in the norm.

Write $v(r,\omega)=v(r\omega)$. The derivative $\partial_rv$ measures variation along a ray. The vector $\nabla_{\Sph}v(r,\omega)$ is the tangential gradient of the function $\omega\mapsto v(r\omega)$ on the unit sphere. It lies in
\[
T_\omega\Sph=\{\xi\in\R^3:\xi\cdot\omega=0\}.
\]
For smooth $v$ it satisfies
\[
\nabla_{\Sph}v(r,\omega)=r(I-\omega\omega^{\mathsf T})\nabla v(r\omega).
\]

The limiting energy is
\begin{align}\label{energy}
\q(v,v)={}&\int_{\B_1}|\nabla v|^2\dd x\\
&+\int_1^R\int_{\Sph}
\Big[\theta r^2|\partial_rv|^2+
 \langle A(\omega)\nabla_{\Sph}v,\nabla_{\Sph}v\rangle\Big]\dd S\dd r.\notag
\end{align}
Here $A(\omega)$ is a symmetric linear map on the two-dimensional plane $T_\omega\Sph$. It specifies the energy of angular variation. Section~\ref{cellappendix} defines it by an explicit minimization over the planar pattern. The estimates needed in the proof are
\begin{equation}\label{coercive}
\frac{\theta}{6}|\xi|^2\leq\langle A(\omega)\xi,\xi\rangle
\leq\theta|\xi|^2,
\qquad \xi\in T_\omega\Sph,
\end{equation}
uniformly in $R$ and $h$. In particular, the lower comparison constant is independent of the small retained fraction $\theta$.

The bilinear energy $\q(v,w)$ is obtained by replacing squares in \eqref{energy} by the corresponding products. The limiting eigenproblem is
\[
\q(v,w)=\lambda\langle v,w\rangle_m
\qquad(w\in H^1(\B_R)).
\]

\textbf{Why angular symmetry is unnecessary.} We do not require $A(\omega)$ to be a scalar multiple of the identity. The weight $m$ and the coefficient of $|\partial_rv|^2$ depend only on $r$, and \eqref{energy} has no product between $\partial_rv$ and $\nabla_{\Sph}v$. These facts let the radial and zero-spherical-average parts be treated separately in Section~\ref{proof}.

\subsection{An explicit definition of the angular matrix}\label{cellappendix}

This definition specifies the coefficient $A$ in \eqref{energy}. The proof in Section~\ref{proof} uses only the two-sided bounds \eqref{coercive}; its identification as the limit is proved with Lemma~\ref{hom}.

Use the planar set $P_{1,h}$ from \eqref{planarset} with $\delta=1$, and put $Y_h=P_{1,h}\cap(0,1)^2$. The notation $H^1_{\mathrm{loc}}(P_{1,h})$ means $H^1$ on every open subset whose closure is a compact subset of $P_{1,h}$. Let $\mathcal V_h$ be the class of functions $\chi\in H^1_{\mathrm{loc}}(P_{1,h})$ such that
\[
\chi(z+a)=\chi(z)\quad(a\in\mathbb Z^2)
\quad\hbox{and}\quad
\int_{Y_h}(|\chi|^2+|\nabla\chi|^2)\dd z<\infty.
\]
The equality expressing periodicity is understood almost everywhere.
Thus the functions in this minimization have period one in the two coordinate directions. All periodicity used below is this explicit equality.

Fix a point $\omega$ in the interior of $K_{i,\sigma}$ and write $\omega=\Psi_{i,\sigma}(y)$. Let
\[
F=D\Psi_{i,\sigma}(y),\qquad G=F^{\mathsf T}F.
\]
Here $F$ is a $3\times2$ matrix of derivatives; its columns span $T_\omega\Sph$. The $2\times2$ matrix $G$ is positive definite. For $\xi\in T_\omega\Sph$, define
\begin{equation}\label{cellformula}
\langle A(\omega)\xi,\xi\rangle
=\inf_{\chi\in\mathcal V_h}
\int_{Y_h}\big|\xi+FG^{-1}\nabla_z\chi(z)\big|^2\dd z.
\end{equation}
The matrices $F,G$ are held fixed at $y$ in this minimization. The formula defines a quadratic form, and hence a symmetric map $A(\omega)$, on the tangent plane. Values on $\Sigma$ may be assigned arbitrarily, since that set has surface area zero.

\textbf{Meaning of the minimization.} The fixed vector $\xi$ represents a prescribed large-scale tangential gradient. The periodic function $\chi$ allows the gradient to vary within one retained planar period. The infimum is the smallest energy obtained from these variations. In the convergence proof, their scaled versions are the small within-period adjustments to the approximating functions.

The choice $\chi=0$ gives $\langle A\xi,\xi\rangle\leq\theta|\xi|^2$. A positive lower bound proportional to $\theta$ follows by integrating along the entire horizontal and vertical strips and using periodicity. The derivatives of the six maps \eqref{psimap} and their inverses are bounded on their closed squares, so the comparison constants are independent of $h$. This is the source of \eqref{coercive}; the proof of Lemma~\ref{hom} establishes this identification.

\subsection{The radial equation and the interface condition}

For a radial function $v(x)=f(|x|)$ define
\begin{equation}\label{p}
p(r)=\begin{cases}r^2,&0<r<1,\\ \theta r^2,&1<r<R.\end{cases}
\end{equation}
Its norm and energy are $4\pi\int_0^Rpf^2$ and $4\pi\int_0^Rp(f')^2$. The equation is
\begin{equation}\label{ode}
-(pf')'=\lambda pf,
\qquad f'(0)=f'(R)=0.
\end{equation}
At $r=1$, $f$ is continuous and
\[
p(1-)f'(1-)=p(1+)f'(1+),
\quad\hbox{equivalently}\quad f'(1-)=\theta f'(1+).
\]
This is the precise meaning of ``flux continuity'': $4\pi pf'$ is the radial derivative multiplied by the available cross-sectional area in the weighted problem. The one-sided weighted derivatives agree at $1$. Both continuity conditions follow from the variational problem on $H^1(\B_R)$.

\section{Auxiliary results}\label{auxiliary}

The six lemmas below provide the geometric, analytic, and radial estimates used in the main proof. Their proofs are given in Section~\ref{lemmproofs}.

\begin{lemma}[Angular geometry]\label{geom}
For the sets in \eqref{Udef}--\eqref{domain}, $U_\delta$ is open, connected, and satisfies $U_\delta=-U_\delta$. For every continuous $g$ on $\Sph$,
\begin{equation}\label{averaging}
\int_{U_\delta}g\dd S\longrightarrow
\theta\int_{\Sph}g\dd S.
\end{equation}
The domains $\Omega_\delta$ are bounded connected Lipschitz domains, starlike about $0$, with
\[
\B_1\subset\Omega_\delta\subset\B_R,
\qquad \partial\Omega_\delta\subset\{x:1\leq|x|\leq R\}.
\]
For fixed $R,\theta$, there are $c,C,a,L>0$, independent of sufficiently small $\delta$, such that:
\begin{itemize}
\item every removed spherical square has diameter between $c\delta$ and $C\delta$, and the distance between distinct removed squares is at least $c\delta$;
\item near every point of $\partial\Omega_\delta$, in a ball of radius $a\delta$ and suitable orthogonal coordinates, the domain is the region on one side of a graph with Lipschitz constant at most $L$.
\end{itemize}
\end{lemma}
\begin{lemma}[Uniform extension]\label{extension}
For fixed $R,\theta$, there are linear operators
\[
E_\delta:H^1(\Omega_\delta)\longrightarrow H^1(\B_R)
\]
such that $E_\delta v=v$ on $\Omega_\delta$ and
\begin{equation}\label{ext}
\|E_\delta v\|_{H^1(\B_R)}\leq C\|v\|_{H^1(\Omega_\delta)},
\end{equation}
with $C$ independent of $\delta$. An extension assigns values on $\B_R\setminus\Omega_\delta$ to a function initially defined on $\Omega_\delta$. The displayed inequality is the required control on those values and their gradients, including near $r=1$ and $r=R$.
\end{lemma}
\noindent The extension estimate has two consequences used below. The extension gives the uniform Sobolev inequality
\[
\|v\|_{L^6(\Omega_\delta)}\leq
C\big(\|\nabla v\|_{L^2(\Omega_\delta)}+\|v\|_{L^2(\Omega_\delta)}\big).
\]
Here and below ``uniform'' means that $C$ does not depend on $\delta$, with $R,\theta$ fixed. Also, the extension must preserve the entire unit ball. Extending only the part in $\mathcal A_{1,R}$ and then attaching it to the ball can give different boundary values at $r=1$. An extension on the whole domain avoids that mismatch. Boundary values of an $H^1$ function, when needed, are understood in the Sobolev sense.

\begin{lemma}[Energy, eigenvalue, and eigenfunction convergence]\label{hom}
For the geometry of Lemma~\ref{geom}, the actual-domain Neumann energies
\[
\q_\delta(v,v)=\int_{\Omega_\delta}|\nabla v|^2\dd x
\]
converge to \eqref{energy}, with mass \eqref{mass}, as follows.
\begin{enumerate}[label=(\roman*)]
\item Every sequence with bounded $L^2$ norm and energy has an extended subsequence converging strongly in $L^2(\B_R)$ to $v\in H^1(\B_R)$. For every such convergent subsequence,
\[
\liminf_\delta\q_\delta(v_\delta,v_\delta)\geq\q(v,v),
\qquad
\int_{\Omega_\delta}v_\delta^2\dd x\longrightarrow\normm{v}^2.
\]
\item For every $v\in H^1(\B_R)$ there is a sequence $v_\delta$ with $E_\delta v_\delta\to v$ in $L^2(\B_R)$ and with convergence of its energy to $\q(v,v)$.
\item The matrix on each tangent plane satisfies \eqref{coercive}; the radial coefficient is exactly $\theta$ and there is no mixed radial--angular term.
\end{enumerate}
Consequently, each indexed Neumann eigenvalue converges to the corresponding limiting eigenvalue, including multiplicities. If the limiting first positive eigenvalue is simple, the normalized first eigenfunctions on the actual domains, after a consistent choice of sign, converge after extension in $L^2(\B_R)$ to the normalized limiting eigenfunction. The actual first positive eigenvalue is then simple for sufficiently small $\delta$.
\end{lemma}
\begin{lemma}[A local boundedness estimate]\label{moser}
For fixed $R,\theta$, there are $C,s_0>0$, independent of $\delta$, with the following property. If $x\in\overline\Omega_\delta$, $0<s\leq s_0$, and $v\in H^1(\Omega_\delta\cap\B(x,2s))$ solves
\[
-\Delta v=F\quad\hbox{in }\Omega_\delta\cap\B(x,2s)
\]
with zero normal derivative on $\partial\Omega_\delta\cap\B(x,2s)$ in the weak sense specified in Section~\ref{proofmoser}, then, for every constant $c\in\R$ and $F\in L^\infty$,
\begin{equation}\label{local}
\|v-c\|_{L^\infty(\Omega_\delta\cap\B(x,s))}
\leq C\Big[s^{-3/2}\|v-c\|_{L^2(\Omega_\delta\cap\B(x,2s))}
 +s^2\|F\|_{L^\infty}\Big].
\end{equation}
The estimate holds also when $x$ is near $r=1$, near a coordinate-region boundary, or near $r=R$. Moreover, for every $\Lambda<\infty$ there is a constant $C=C(R,\theta,\Lambda)$ such that every Neumann eigenfunction normalized by $\int_{\Omega_\delta}u_\delta^2\,dx=1$, with eigenvalue at most $\Lambda$, satisfies
\begin{equation}\label{globalbound}
\|u_\delta\|_{L^\infty(\Omega_\delta)}\leq C.
\end{equation}
Each has a continuous representative on $\overline\Omega_\delta$; the same local estimate holds for its boundary values.
\end{lemma}
\begin{lemma}[A small radial eigenvalue]\label{trial}
Let $\lambda_{\rm rad}$ denote the first positive eigenvalue of \eqref{energy} restricted to radial functions. With \eqref{theta},
\begin{equation}\label{radbound}
\lambda_{\rm rad}\leq C R^{-5/2}
\end{equation}
for large $R$, with $C$ independent of $R$.

\end{lemma}
\begin{lemma}[The first radial eigenfunction]\label{sturm}
Put
\[
p(r)=r^2\ (0<r<1),\qquad p(r)=\theta r^2\ (1<r<R).
\]
The first positive radial eigenvalue is simple. A corresponding real eigenfunction $f$, chosen with $f(0)>0$, has exactly one zero in $(0,R)$, has weighted mean zero, and satisfies \eqref{ode} with equal one-sided values of $pf'$ at $r=1$. In the unit ball,
\begin{equation}\label{sinc}
f(r)=f(0)\frac{\sin(\sqrt\lambda\,r)}{\sqrt\lambda\,r},
\end{equation}
with the value at $r=0$ understood by continuity.
\end{lemma}
\section{Proof of the main theorem}\label{proof}

We prove more precisely that, for all sufficiently large fixed $R$ and then all sufficiently small $\delta$, the domains $\Omega_\delta$ defined in \eqref{domain} satisfy every conclusion of Theorem~\ref{main}.

\subsection{The first limiting eigenfunction is radial}

For $v\in H^1(\B_R)$, define its spherical average
\[
\bar v(r)=\frac1{4\pi}\int_{\Sph}v(r,\omega)\dd S,
\qquad v^\perp=v-\bar v.
\]
Spherical averaging is a bounded projection on $H^1(\B_R)$: Jensen's inequality controls its $L^2$ norm and radial derivative, and its angular derivative vanishes. It also preserves the weighted mean. Thus $v^\perp$ has zero spherical average for almost every $r$. Radial functions and these mean-zero functions are orthogonal for both the mass inner product and the energy. For the radial energy term, the cross term vanishes because $\int_{\Sph}\partial_rv^\perp\,dS=0$; for the angular term, it vanishes because $\nabla_{\Sph}\bar v=0$. The same statements hold in the ball.

Consequently eigenvalue minimization can be carried out separately among radial functions and among functions with zero spherical average, even when $A(\omega)$ depends on the angular direction. For the latter functions, the spherical Poincar\'e inequality is
\[
\int_{\Sph}|\nabla_{\Sph}v^\perp|^2\dd S
\geq2\int_{\Sph}|v^\perp|^2\dd S.
\]
For $0<r<1$, $r^2\leq1$; for $1<r<R$, $r^2\leq R^2$. With \eqref{coercive}, this gives
\begin{equation}\label{angularbound}
\q(v^\perp,v^\perp)\geq
\frac{1}{3R^2}\normm{v^\perp}^2.
\end{equation}
Lemma~\ref{trial} gives $\lambda_{\rm rad}\leq CR^{-5/2}$. Choose $R$ so large that
\[
CR^{-5/2}<\frac{1}{3R^2}.
\]
The first positive limiting eigenvalue is therefore $\lambda=\lambda_{\rm rad}$. It is simple by Lemma~\ref{sturm}, and its eigenfunction $u_0(x)=f(|x|)$ is radial. Normalize it by $\normm{u_0}=1$ and $f(0)>0$.

\subsection{Strict radial decrease}

Let $z\in(0,R)$ be its unique zero. Then $f>0$ before $z$ and $f<0$ after $z$. Define
\[
J(r)=\int_0^r p(s)f(s)\dd s.
\]
For $0<r\leq z$, $J(r)>0$. For $z<r<R$, the weighted mean-zero property gives
\[
J(r)=-\int_r^R p(s)f(s)\dd s>0.
\]
Integrating \eqref{ode} from $0$ to $r$, using $pf'\to0$ at $0$ and equality of the one-sided weighted derivatives at $1$, gives
\begin{equation}\label{decrease}
p(r)f'(r)=-\lambda J(r)<0\qquad(0<r<R).
\end{equation}
Thus $f$ is strictly decreasing, with one-sided derivatives at $1$. In particular,
\begin{equation}\label{gap}
G=f(0)-f(1)>0,
\qquad \sup_{1\leq r\leq R}f(r)=f(1).
\end{equation}
Here $R$ is fixed, so the positive number $G$ is independent of $\delta$.

\subsection{Eigenvalue and eigenfunction convergence}

By Lemma~\ref{hom}, for sufficiently small $\delta$ the first positive actual-domain eigenvalue $\lambda_\delta$ is simple and tends to $\lambda$. Choose its real eigenfunction $u_\delta$ with
\[
\int_{\Omega_\delta}u_\delta\dd x=0,
\qquad \int_{\Omega_\delta}u_\delta^2\dd x=1,
\]
and choose its sign to converge to $u_0$. Put
\[
\eta_\delta=\|E_\delta u_\delta-u_0\|_{L^2(\B_R)}\longrightarrow0.
\]
The eigenvalues are bounded and Lemma~\ref{moser} supplies a uniform bound on $\|u_\delta\|_\infty$. Hence
\[
F_\delta=\lambda_\delta u_\delta
\]
is bounded in $L^\infty$, independently of $\delta$.

\subsection{Uniform convergence on the actual domains}

Define
\[
\omega_0(t)=\sup\{|u_0(x)-u_0(y)|:
x,y\in\overline\B_R,\ |x-y|\leq t\}.
\]
Since $u_0$ is continuous on the compact ball, $\omega_0(t)\to0$ as $t\downarrow0$. Fix $s\in(0,s_0]$.

For any $x\in\overline\Omega_\delta$, apply \eqref{local} to $v=u_\delta$ and the \emph{constant} $c=u_0(x)$. In $\Omega_\delta\cap\B(x,2s)$,
\begin{align*}
\|u_\delta-u_0(x)\|_{L^2}
&\leq \|u_\delta-u_0\|_{L^2}
 +\|u_0-u_0(x)\|_{L^2}\\
&\leq\eta_\delta+C s^{3/2}\omega_0(2s).
\end{align*}
The local estimate therefore gives, with a new constant $C$,
\begin{equation}\label{uniformestimate}
\sup_{x\in\overline\Omega_\delta}|u_\delta(x)-u_0(x)|
\leq C\Big[s^{-3/2}\eta_\delta+\omega_0(2s)+s^2\Big].
\end{equation}
First let $\delta\downarrow0$ with $s$ fixed; then let $s\downarrow0$. It follows that
\begin{equation}\label{uniformconv}
\sup_{x\in\overline\Omega_\delta}|u_\delta(x)-u_0(x)|\longrightarrow0.
\end{equation}

The use of a constant preserves the zero normal derivative: $u_\delta-u_0(x)$ has the same boundary condition as $u_\delta$. The estimate is applied on the actual domain, including its boundary.

\subsection{The strict boundary gap}

Take $\delta$ so small that the error in \eqref{uniformconv} is less than $G/4$. Lemma~\ref{geom} and \eqref{gap} give
\[
u_0(x)\leq f(1)=f(0)-G\qquad(x\in\partial\Omega_\delta).
\]
Therefore
\[
u_\delta(0)\geq f(0)-G/4,
\qquad
\max_{\partial\Omega_\delta}u_\delta\leq f(0)-3G/4.
\]
This proves \eqref{target}, with a gap at least $G/2$.

\subsection{The unique maximum is at the origin}

Because $\Omega_\delta=-\Omega_\delta$, the function $u_\delta(-x)$ is another eigenfunction for the same simple eigenvalue. Thus it equals either $u_\delta(x)$ or $-u_\delta(x)$. For small $\delta$, \eqref{uniformconv} gives $u_\delta(0)>0$, which excludes the minus sign. Hence $u_\delta$ is even and
\[
\nabla u_\delta(0)=0.
\]

Put $w_\delta=u_\delta-u_0$. On the fixed ball $\B_1$,
\[
-\Delta w_\delta
=\lambda_\delta w_\delta+(\lambda_\delta-\lambda)u_0.
\]
Choose fixed $0<\rho<\rho_1<\rho_2<1$. The interior $H^2$ estimate for the Laplacian gives
\[
\|w_\delta\|_{H^2(\B_{\rho_1})}
\leq C\big(\|w_\delta\|_{L^2(\B_{\rho_2})}
             +|\lambda_\delta-\lambda|\big)\longrightarrow0.
\]
The corresponding estimate for two more derivatives, on the smaller ball, gives
\[
\|w_\delta\|_{H^4(\B_\rho)}
\leq C\big(\|w_\delta\|_{H^2(\B_{\rho_1})}
             +|\lambda_\delta-\lambda|\big)\longrightarrow0.
\]
Here $u_0$ is smooth in $\B_1$ and the eigenvalues are bounded. In dimension three, $H^4(\B_\rho)$ embeds continuously into $C^2(\overline\B_\rho)$, so
\[
u_\delta\longrightarrow u_0\quad\hbox{in }C^2(\overline\B_\rho)
\qquad(\rho<1).
\]

Formula~\eqref{sinc} gives
\[
D^2u_0(0)=-\frac{\lambda f(0)}3 I.
\]
Choose $0<\rho<1$ and $a>0$ so that $D^2u_0\leq-2aI$ throughout $\overline\B_\rho$. For sufficiently small $\delta$, $D^2u_\delta\leq-aI$ there. Since $\nabla u_\delta(0)=0$,
\[
u_\delta(x)\leq u_\delta(0)-\frac a2|x|^2
\qquad(0<|x|\leq\rho).
\]
On the rest of the domain, strict radial decrease gives the fixed gap
\[
H=f(0)-f(\rho)>0.
\]
Take the uniform error in \eqref{uniformconv} smaller than $H/3$. For $|x|\geq\rho$,
\[
u_\delta(x)\leq f(\rho)+H/3
=f(0)-2H/3<u_\delta(0).
\]
Combining the two regions proves \eqref{exact} and completes the proof.\hfill$\square$

\section{Proofs of the auxiliary results}\label{lemmproofs}

We use standard extension, trace, Sobolev, compactness, and Poincar\'e results on fixed bounded Lipschitz domains; see \cite{Weber,Evans}. The constructions below establish the estimates uniform in $\delta$ that are needed here. We also use the variational characterization of Neumann eigenvalues, the spherical Poincar\'e inequality with constant $2$, and interior elliptic estimates. Boundary continuity follows from \cite[Proposition~3.6]{Nittka}, as used in Section~\ref{proofmoser}.

\subsection{Proof of Lemma~\ref{geom}}\label{proofgeom}

Throughout this proof $R$ and $h$ are fixed. All constants may depend on
them, but not on $\delta$. Put $Q=(0,1)^2$ and
\[
H=[h/2,1-h/2]^2,\qquad Y=Q\setminus H,
\qquad Q_a=\delta(a+Q).
\]
Thus $H_a=\delta(a+H)$. The six maps $\Psi_{i,\sigma}$ and their
inverses have uniformly bounded first derivatives on their closed
coordinate squares. Consequently there are absolute constants
$b_0,b_1>0$ such that
\begin{equation}\label{geom-bilip}
b_0|y-z|\leq d_{\Sph}(\Psi_{i,\sigma}(y),\Psi_{i,\sigma}(z))
\leq b_1|y-z|\qquad(y,z\in[-1,1]^2).
\end{equation}
For the upper bound one can take $b_1=1$: the norm of $D\Psi$ is at
most $1$, and the image of the straight segment from $y$ to $z$ has
length at most $|y-z|$. The lower bound follows, for example, from the
displayed inverse formula and $|\omega_i|\geq1/\sqrt3$ on the face.

Every selected $H_a$ is at spherical distance greater than $4\delta$
from $\Sigma$. Its surrounding entire period square $\overline Q_a$
is still inside the coordinate face. Indeed, every point of
$\overline Q_a$ is at Euclidean distance at most $h\delta/\sqrt2$
from $H_a$. If $\overline Q_a$ crossed a face boundary, the line
segment from $H_a$ to the first boundary point, followed by $\Psi$,
would contradict the $4\delta$ distance condition. In particular,
\begin{equation}\label{period-far-from-seam}
d_{\Sph}(\Psi(\overline Q_a),\Sigma)>3\delta.
\end{equation}
The images of these surrounding period squares have disjoint
interiors on each face, and they are disjoint across different faces.
This additional observation will be useful for the extension proof.

The removed sets are finitely many disjoint closed topological disks.
Their diameters lie between $b_0(1-h)\delta$ and
$b_1\sqrt2(1-h)\delta$. On the same face, distinct squares have
coordinate distance at least $h\delta$, so their spherical distance
is at least $b_0h\delta$. For squares on distinct faces, every path
from one to the other crosses $\Sigma$, and their distance is
greater than $8\delta$. Thus the asserted separation estimate holds.
Their complement is open and connected. To verify connectedness
directly, join two points of the complement by a path on the sphere.
Enlarge each removed square slightly to a closed concentric square
$H_a^+$, still within its face, such that the enlarged images are
mutually disjoint and contain neither endpoint. This is possible
by separation and the positive distance of each endpoint from the
finitely many removed sets. Whenever the path enters the image of
$H_a^+$, replace the portion between its first and last intersections
with that image by a path along $\Psi(\partial H_a^+)$. This curve
lies entirely in $U_\delta$. The enlarged images are disjoint, so
a replacement cannot meet another removed set. There are only
finitely many disks, and these replacements give a path in
$U_\delta$.
Antipodal symmetry follows from the two identities for $\Psi$ and
$H_a$ already displayed in the construction.

Here is the averaging argument, including the effect of the removal
rule. In each coordinate square the ordinary repeated pattern has
indicator $\mathbf1_Y(y/\delta)$, interpreted periodically. For every
continuous $f$ on $[-1,1]^2$,
\begin{equation}\label{periodic-average-proof}
\int_{[-1,1]^2}f(y)\mathbf1_Y(y/\delta)\,dy
\longrightarrow \theta\int_{[-1,1]^2}f(y)\,dy.
\end{equation}
To prove this, sum over the whole period squares contained in
$[-1,1]^2$. On each such square subtract the value of $f$ at one
point. The error is at most its area times the modulus of continuity
of $f$ at $\sqrt2\delta$. The squares intersecting the boundary have
total area $O(\delta)$. The integral of the periodic indicator on a
whole square is exactly $\theta\delta^2$, which proves
\eqref{periodic-average-proof}.

A complementary square not removed by the rule either intersects
the boundary of its face or has an image meeting
$\{d_{\Sph}(\omega,\Sigma)\leq4\delta\}$. Its whole image is
therefore within distance $(4+\sqrt2)\delta$ of $\Sigma$, by
\eqref{geom-bilip}. The face-boundary squares satisfy the same type
of bound. Since $\Sigma$ is a finite union of smooth arcs, this
neighborhood has area $O(\delta)$: an arc of length $\ell$ has a
neighborhood of width $t$ with area at most $C(\ell t+t^2)$, as is
seen by covering it with consecutive geodesic rectangles and two
end disks. In \eqref{periodic-average-proof} take
\[
f(y)=g(\Psi(y))J(y),\qquad
J(y)=(1+y_1^2+y_2^2)^{-3/2},
\]
and sum over the six faces. The retained extra neighborhood has
vanishing area. This proves \eqref{averaging}.

The domain formula immediately shows that $\Omega_\delta$ is open,
bounded, starlike, and contains $\B_1$. Any point can be joined to
$0$ by its radial segment, so the domain is connected. Its boundary
can contain no point with $|x|<1$ or $|x|>R$, proving the asserted
boundary inclusion.

For clarity, we give the local boundary check rather than inferring
it from openness. Inside one face, a removed set is bounded by four
planes through $0$: an edge $y_1=c$ becomes
$x_j=c\sigma x_i$, and an edge $y_2=c$ becomes
$x_k=c\sigma x_i$. Its two radial ends are the spheres $r=1$ and
$r=R$. Near a corner at the lower end, suitable signed coordinates
$(s,z_1,z_2)$ make the removed set
\[
\{s\geq0,\ z_1\geq0,\ z_2\geq0\}.
\]
Here $s=r-1$; one or both $z$ inequalities are omitted at an edge
or in the relative interior of a face. The domain is the complement
of this orthant. At the upper end, with $s=R-r$, the corresponding
domain is
\[
\{s>0\}\cap\big(\{z_1<0\}\cup\{z_2<0\}\big).
\]
Again one inequality is omitted when only one edge is present.
Away from the ends the same models have an unrestricted $s$.

All these models are Lipschitz graph domains. More explicitly, in
coordinates $q=t(1,1,1)/\sqrt3+z$ with $z_1+z_2+z_3=0$, the
complement of the positive orthant is
\[
t<\max_i(-\sqrt3 z_i).
\]
The upper-end model, after reversing the signs of the last two
coordinates, is an intersection of one positive half-space and
the union of two positive half-spaces; in the same diagonal coordinates its graph is explicitly
\[
t>\max\big(-\sqrt3 z_1,\min(-\sqrt3 z_2,-\sqrt3 z_3)\big).
\]
Maxima and minima preserve the Lipschitz bound of these affine functions. The two-dimensional
edge models are the corresponding formulas with two coordinates.

The actual coordinates $(r,y_1,y_2)\mapsto r\Psi(y)$ are smooth
with bounded invertible derivative for $1/2\leq r\leq R+1$.
On a neighborhood of diameter $O(\delta)$ their derivative differs
from its value at the center by $O(\delta)$, uniformly. A linear
invertible image of any of the preceding graph models is a graph
in the direction which is the image of its original graph
direction, with a bound depending only on the derivative and its
inverse. The $O(\delta)$ perturbation retains this property for
sufficiently small $\delta$. To check the last assertion, project
the perturbed boundary onto the plane perpendicular to this
direction. The projection is the sum of an invertible linear map
and a map with Lipschitz norm $O(\delta)(1+L_0)$, where $L_0$
is the model graph bound. Its inverse and the height function are
therefore Lipschitz with uniform bounds.

Choose $a>0$ smaller than a fixed multiple of both $h$ and $1-h$.
The identity
$|r\omega-s\eta|^2=(r-s)^2+rs|\omega-\eta|^2$, for $r,s\geq1$,
transfers the angular separation estimate to the excluded radial
sets. Thus a ball of radius $a\delta$ meets at most one removed
square, and cannot meet two opposite
edges of that square. Choose also $\delta$ small enough that it
cannot meet both radial ends. It thus has one of the graph models
just listed. Near $\Sigma$ no removed square is present by
\eqref{period-far-from-seam}; only the smooth outer sphere can be a
boundary there. This proves the quantified Lipschitz statement and
completes the geometry lemma.\hfill$\square$

\subsection{Proof of Lemma~\ref{extension}}\label{proofextension}

We construct the operator. An angular extension used separately
above $r=1$ would not preserve the trace of the function on the
unit sphere. The following construction includes a three-dimensional
extension on a layer of thickness $\delta$ at this interface.

We first record precisely the two fixed-domain operators used.
Besides $Q,H,Y$ above, define the open three-dimensional set
\begin{equation}\label{bottom-model-D}
D=\big((-1,0)\times Q\big)\ \cup\
\big((-1,2)\times Y\big),\qquad F=(-1,2)\times Q.
\end{equation}
There are linear operators $T:H^1(Y)\to H^1(Q)$ and
$S:H^1(D)\to H^1(F)$, equal to their input on its original set,
preserving constants, and satisfying
\begin{align}\label{fixed-T-estimates}
\|Tu\|_{L^2(Q)}&\leq C_h\|u\|_{L^2(Y)},&
\|\nabla Tu\|_{L^2(Q)}&\leq C_h\|\nabla u\|_{L^2(Y)},\\
\label{fixed-S-estimates}
\|Su\|_{L^2(F)}&\leq C_h\|u\|_{L^2(D)},&
\|\nabla Su\|_{L^2(F)}&\leq C_h\|\nabla u\|_{L^2(D)}.
\end{align}

Both $Y$ and $D$ are fixed bounded connected Lipschitz domains; their
corner geometry was checked above. The fixed-domain extension theorem
and its reflection construction \cite[Theorem~5.2.1 and its proof]{Weber}
give, for either source $G$, a linear extension $A$ into the surrounding
rectangle with
\begin{equation}\label{raw-reflection-estimates}
\|Au\|_2\leq C_h\|u\|_2,
\qquad \|\nabla Au\|_2\leq C_h(\|\nabla u\|_2+\|u\|_2).
\end{equation}
The separate $L^2$ bound follows from the same reflection and partition
of unity, without differentiating. If $u_G$ is the mean of $u$ on $G$,
the Poincar\'e inequality \cite[Theorem~7.1.1]{Weber} gives
$\|u-u_G\|_2\leq C_h\|\nabla u\|_2$.
Thus replacing $A$ by $u\mapsto A(u-u_G)+u_G$ preserves the extension
property and the $L^2$ bound, preserves constants, and removes the
$\|u\|_2$ term from the gradient bound. This proves
\eqref{fixed-T-estimates}--\eqref{fixed-S-estimates}.
In particular, $T$ is bounded on $L^2$ and hence commutes with weak
radial differentiation, by difference quotients:
\begin{equation}\label{T-radial-commutation}
\partial_r(Tu)=T(\partial_r u).
\end{equation}

We now apply these fixed operators to each selected square. Fix one
such square on one face and suppress its face index. Use coordinates
\[
x=r\Psi(y),\qquad y=\delta(a+z),\quad z\in Q.
\]
For $r\geq1$, the source over this period is precisely $z\in Y$.
For $r<1$ the whole period is in the domain. Write
$\hat v(r,y)=v(r\Psi(y))$. Let $b_a(r,y)$ be the scaled angular
extension $T$ of $\hat v(r,\cdot)$, for $1<r<R$.
For $1-\delta<r<1+2\delta$, let $s_a(r,y)$ be the scaled
three-dimensional extension $S$ under
\[
t=(r-1)/\delta,\qquad z=y/\delta-a.
\]
Its source is exactly the image of $D$ in \eqref{bottom-model-D}.
For all sufficiently small $\delta$ this image lies inside $\B_R$.
By \eqref{period-far-from-seam}, it meets no other removed set.

Fix a smooth function $\chi$ on $\R$ with $0\leq\chi\leq1$,
$\chi(t)=1$ for $t\leq1$, and $\chi(t)=0$ for $t\geq3/2$.
The function $\chi$ is a cutoff: in the formula below it selects $s_a$ for $r\leq1+\delta$ and $b_a$ for $r\geq1+3\delta/2$. On the excluded set associated with $H_a$, define
\begin{equation}\label{extension-formula}
(E_\delta v)(r\Psi(y))=
\begin{cases}
\chi((r-1)/\delta)s_a(r,y)
 +(1-\chi((r-1)/\delta))b_a(r,y),&1\leq r<1+2\delta,\\
b_a(r,y),&1+2\delta\leq r<R.
\end{cases}
\end{equation}
On $\Omega_\delta$ itself set $E_\delta v=v$.
Because both local operators equal the input on $Y$, the formulas
have the original trace on every lateral edge. At $r=1$ the cutoff
is $1$, and $S$ equals $v$ throughout the lower half of its box;
hence it gives exactly the unit-ball trace. In fact near the outer
edges of each period square the result equals $v$ on an open strip,
since $H_a$ is inset by $h\delta/2$. Thus different period
definitions glue without a jump. At $r=R$ no further gluing is
needed: the target space is $H^1(\B_R)$. These observations and
the trace gluing rule for $H^1$ functions prove that the result
belongs to $H^1(\B_R)$. The entire unit ball is preserved.

We give the uniform estimates, including the term created by the
cutoff. In these coordinates
\[
dx=r^2J(y)\,dr\,dy,
\]
and the gradient norm is uniformly equivalent, for fixed $R$ and
$1/2<r<R$, to
\[
\int (|\partial_r\hat v|^2+|\nabla_y\hat v|^2)\,dr\,dy.
\]
All equivalence constants are independent of $\delta$ and $a$.
Scaling \eqref{fixed-T-estimates} and using
\eqref{T-radial-commutation} give
\begin{align*}
\int_{Q_a}|b_a|^2\,dy&\leq C_h\int_{Q_a\setminus H_a}|\hat v|^2\,dy,\\
\int_{Q_a}|\nabla_y b_a|^2\,dy&\leq C_h\int_{Q_a\setminus H_a}|\nabla_y\hat v|^2\,dy,\\
\int_{Q_a}|\partial_r b_a|^2\,dy&\leq C_h\int_{Q_a\setminus H_a}|\partial_r\hat v|^2\,dy.
\end{align*}
The scaled $S$ has the same $L^2$ and full-gradient bounds on
its three-dimensional box, because all three coordinates are
scaled by $\delta$.

Let $D_a$ denote its source box in $(r,y)$ coordinates, and let
$c_a$ be the ordinary coordinate mean of $\hat v$ on $D_a$.
Its scaled Poincar\'e inequality reads
\[
\int_{D_a}|\hat v-c_a|^2\,dr\,dy
\leq C_h\delta^2\int_{D_a}
(|\partial_r\hat v|^2+|\nabla_y\hat v|^2)\,dr\,dy.
\]
Since $S$ and $T$ preserve constants and are $L^2$ bounded, on
the transition box $1+\delta<r<1+3\delta/2$ this implies
\begin{equation}\label{blend-difference-estimate}
\int |s_a-b_a|^2\,dr\,dy
\leq C_h\delta^2\int_{D_a}
(|\partial_r\hat v|^2+|\nabla_y\hat v|^2)\,dr\,dy.
\end{equation}
The extra derivative in \eqref{extension-formula} is
$\delta^{-1}\chi'((r-1)/\delta)(s_a-b_a)$.
Equation \eqref{blend-difference-estimate} bounds its squared
integral by the source energy, with no negative power of
$\delta$ remaining.

The angular period interiors are disjoint, including across faces.
The source sets $D_a$ are therefore disjoint as well. Within each
period the angular and bottom source estimates overlap at most
twice. Summing the estimates and adding the unchanged part of the
domain proves the stronger pair of bounds
\begin{equation}\label{extension-strong-pair}
\|E_\delta v\|_{L^2(\B_R)}\leq C_{R,h}\|v\|_{L^2(\Omega_\delta)},
\qquad
\|\nabla E_\delta v\|_{L^2(\B_R)}
\leq C_{R,h}\|\nabla v\|_{L^2(\Omega_\delta)}.
\end{equation}
This proves \eqref{ext}. Applying the ordinary Sobolev inequality
on the fixed ball to $E_\delta v$ gives the stated uniform
$L^6$ inequality on $\Omega_\delta$.

\hfill$\square$

\subsection{Proof of Lemma~\ref{hom}}
\label{proofhom}

Fix $R$ and $h$. Compactness and approximation constants may depend on them; the bounds on $A/\theta$ below do not.

\paragraph{1. Coordinates and the exact energy formula.}
Put $Q=(-1,1)^2$, $Y=(0,1)^2$, and
\[
Y_h=P_{1,h}\cap Y,\qquad \chi_h=\mathbf1_{P_{1,h}}.
\]
The function $\chi_h$ is regarded as a function of period $1$ in each coordinate. Its average over $Y$ is $\theta$.
On one of the six spherical regions, suppress the two face indices and write
\[
F(y)=D\Psi(y),\quad G(y)=F(y)^{\mathsf T}F(y),\quad
J(y)=\sqrt{\det G(y)},\quad X(r,y)=r\Psi(y).
\]
The functions $G,G^{-1},J,J^{-1}$ are smooth and bounded on $\overline Q$, uniformly over the six faces. For $w=v\circ X$, orthogonality of $\Psi$ to the columns of $F$ gives
\begin{equation}\label{h-exact}
dx=r^2J(y)\,dr\,dy,\qquad
|\nabla v(X(r,y))|^2
=|\partial_rw|^2+r^{-2}(\nabla_yw)^{\mathsf T}G(y)^{-1}\nabla_yw.
\end{equation}
In particular, $\nabla_{\Sph}v=FG^{-1}\nabla_yw$.

For every compact subset $Q'\subset Q$, the retained set in these coordinates is, for all sufficiently small $\delta$, exactly
\begin{equation}\label{h-local-pattern}
\mathbf1_{U_\delta}(\Psi(y))=\chi_h(y/\delta)\qquad(y\in Q').
\end{equation}
Indeed, every complementary square meeting $Q'$ eventually satisfies the removal rule. The additionally retained part lies within distance $C\delta$ of $\partial Q$ and has spherical area $O(\delta)$. Periodic averaging, first for continuous tests and then by $L^1$ approximation, therefore gives
\begin{equation}\label{h-weight}
\mathbf1_{\Omega_\delta}\ \stackrel{*}{\rightharpoonup}\ m
\quad\hbox{in }L^\infty(\B_R).
\end{equation}

\paragraph{2. The cell minimization and its bounds.}
For a positive definite $2\times2$ matrix $D$ and $\eta\in\R^2$, put
\begin{equation}\label{h-cell-D}
B(D,\eta)=\min_{\zeta\in\mathcal V_h}
\int_{Y_h}(\eta+\nabla_z\zeta)^{\mathsf T}D(\eta+\nabla_z\zeta)\,dz.
\end{equation}
On the connected Lipschitz set $Y_h$, impose equal traces on opposite outside edges and $\int_{Y_h}\zeta=0$. These define a closed subspace of $H^1(Y_h)$. Poincar\'e's inequality and the Lax--Milgram theorem \cite{Evans} give a unique minimizer, characterized by
\begin{equation}\label{h-cell-weak}
\int_{Y_h}\nabla\varphi^{\mathsf T}D(\eta+\nabla\zeta)\,dz=0
\qquad(\varphi\in\mathcal V_h).
\end{equation}
The minimizer depends linearly on $\eta$. Filling the central square by a bounded $H^1$ extension leaves the outside-edge traces unchanged, and hence extends cell functions to $H^1_{\rm per}(Y)$. Periodic convolution then gives density of smooth periodic restrictions in $\mathcal V_h$. This fixed-cell extension may depend on $h$.

If $d_-I\leq D\leq d_+I$, the following elementary estimates hold:
\begin{equation}\label{h-cell-bounds}
d_-h|\eta|^2\leq B(D,\eta)\leq\theta\eta^{\mathsf T}D\eta.
\end{equation}
For the lower bound, fix $z_2$ in the horizontal strip, of total width $h$. Periodicity gives
\[
\int_0^1(\eta_1+\partial_{z_1}\zeta)\,dz_1=\eta_1,
\]
so Cauchy--Schwarz gives an integral of its square at least $\eta_1^2$. Integrating in this strip gives $h\eta_1^2$. The vertical strip gives $h\eta_2^2$ for the other component. These two estimates can be added because they concern the two different components of $|\eta+\nabla\zeta|^2$, even in their overlap. Multiply by $d_-$ and take the infimum. The upper bound uses $\zeta=0$. The argument holds for all $H^1$ cell functions by Fubini and their periodic trace values.

For a tangent vector $\xi$ at $\omega=\Psi(y)$ put $\eta=F(y)^{\mathsf T}\xi$. Then $FG^{-1}\eta=\xi$, and the definition in Section~\ref{cellappendix} is precisely
\begin{equation}\label{h-cell-A}
\xi^{\mathsf T}A(\omega)\xi=B(G(y)^{-1},\eta).
\end{equation}
For completeness, this gives the explicit constants in \eqref{coercive}. Put $s=|y|^2$. Direct differentiation of \eqref{psimap} gives
\[
G=(1+s)^{-1}I-(1+s)^{-2}yy^{\mathsf T}.
\]
Its eigenvalues are $(1+s)^{-1}$ and $(1+s)^{-2}$, so those of $G^{-1}$ are $1+s$ and $(1+s)^2$. Using \eqref{h-cell-bounds} and $|\xi|^2=\eta^{\mathsf T}G^{-1}\eta$, we obtain
\[
\xi^{\mathsf T}A(\omega)\xi
\geq h(1+s)|\eta|^2
\geq\frac{h}{1+s}|\xi|^2
\geq\frac{\theta}{6}|\xi|^2.
\]
Here $s\leq2$ and $h\geq\theta/2$. The upper bound is $\xi^{\mathsf T}A(\omega)\xi\leq\theta|\xi|^2$, by the zero cell function. Both constants are independent of $h$ and $R$.

Write $\zeta(D,\eta)$ for the normalized minimizer. Subtracting two instances of \eqref{h-cell-weak} and testing with the difference of the minimizers gives, on compact sets of positive definite matrices,
\[
\|\nabla\zeta(D,\eta)-\nabla\zeta(\widetilde D,\widetilde\eta)\|_2
\leq C\big(|\eta-\widetilde\eta|+|D-\widetilde D|\,|\widetilde\eta|\big).
\]
By Poincar\'e, this also gives $H^1$ continuity. Thus $A$ is continuous within each face and measurable on the sphere; its values on the face boundaries are irrelevant.

\paragraph{3. Compactness and convergence of the norm.}
For $v_\delta$ with bounded norm and energy, Lemma~\ref{extension} and Rellich compactness give, along a subsequence, $\widetilde v_\delta=E_\delta v_\delta\rightharpoonup v$ in $H^1(\B_R)$ and $\widetilde v_\delta\to v$ in $L^2(\B_R)$. Hence $\widetilde v_\delta^2\to v^2$ in $L^1$, and \eqref{h-weight} gives
\begin{equation}\label{h-mass-conv}
\int_{\Omega_\delta}v_\delta^2\,dx\longrightarrow\int_{\B_R}mv^2\,dx.
\end{equation}
Polarization gives convergence of the corresponding inner products. Recall that $m\geq\theta>0$.

\paragraph{4. What the gradients can do within a period.}
We give the needed partial unfolding argument; see \cite{CDG} for the general method. Here $H^1_{\rm per}(Y)$ consists of $H^1(Y)$ functions with equal opposite-edge traces.

Suppose $w_\delta$ is bounded in $H^1(I\times Q)$, where $I$ is a bounded interval, and $w_\delta\to w$ in $L^2$. For every square $\delta(a+Y)$ contained in $Q$, and every $y$ in that square, define
\[
T_\delta f(r,y,z)=f(r,\delta a+\delta z),\qquad z\in Y.
\]
Set $T_\delta f=0$ outside the union $Q_\delta$ of these squares. Change of variables gives
\[
\int_{I\times Q\times Y}|T_\delta f|^2\,dr\,dy\,dz
=\int_{I\times Q_\delta}|f|^2\,dr\,dy.
\]
Thus $T_\delta$ is an $L^2$ contraction. By continuous approximation, $T_\delta w\to w$, so also $T_\delta w_\delta\to w$ strongly in $L^2$. Extract weak $L^2$ limits $q_r,q_y$ of the unfolded derivatives.

For smooth compactly supported $\varphi(r,y)$ and smooth periodic $\psi(z)$,
\begin{equation}\label{h-test-rel}
\int_{I\times Q}f_\delta(r,y)\varphi(r,y)\psi(y/\delta)\,dr\,dy
\longrightarrow
\int_{I\times Q\times Y}q(r,y,z)\varphi(r,y)\psi(z)\,dr\,dy\,dz
\end{equation}
whenever $f_\delta$ is bounded in $L^2$ and $T_\delta f_\delta\rightharpoonup q$. Indeed, freezing $\varphi$ within each square introduces an error at most $C\delta\|f_\delta\|_2$; the omitted boundary region contributes at most $C|Q\setminus Q_\delta|^{1/2}\|f_\delta\|_2$. This also proves periodic averaging, with integrable periodic factors handled by $L^1(Y)$ approximation.

Integration by parts gives
\[
\int\partial_rw_\delta\,\varphi\,\psi(y/\delta)
=-\int w_\delta\,\partial_r\varphi\,\psi(y/\delta).
\]
Strong $L^2$ convergence, periodic averaging, and \eqref{h-test-rel} yield
\begin{equation}\label{h-radial-limit}
q_r(r,y,z)=\partial_rw(r,y).
\end{equation}
For the angular derivatives, constant $\psi$ first gives
\[
\int_Y q_y(r,y,z)\,dz=\nabla_yw(r,y).
\]
For periodic vector fields with $\operatorname{div}_z\psi=0$, integration by parts has no $\delta^{-1}$ term:
\[
\int\nabla_yw_\delta\cdot\varphi\psi(y/\delta)
=-\int w_\delta\,\nabla_y\varphi\cdot\psi(y/\delta).
\]
Passing to the limit shows that $b=q_y-\nabla_yw$ is orthogonal in the periodic variable to every divergence-free field. Its Fourier coefficients satisfy $\widehat b(0)=0$ and $\widehat b(k)\parallel k$ for $k\ne0$. Define $W$ by
\[
\widehat W(k)=\frac{k\cdot\widehat{(q_y-\nabla_yw)}(k)}{2\pi i|k|^2},
\quad(k\ne0),\qquad \widehat W(0)=0.
\]
Parseval's identity gives $W\in L^2(I\times Q;H^1_{\rm per}(Y))$ and
\begin{equation}\label{h-angular-limit}
q_y=\nabla_yw+\nabla_zW.
\end{equation}
A countable family of Fourier tests establishes these identities for almost every $(r,y)$; the formula also gives measurability of $W$.

\paragraph{5. The lower energy bound.}
Compose the extended functions $\widetilde v_\delta$ with $X(r,y)$ on each face. Formula \eqref{h-exact} shows that the resulting functions are bounded in $H^1((1,R)\times Q)$ and converge strongly in $L^2$ to $w=v\circ X$. Fix a smaller rectangle $Q'\subset Q$ and a smaller interval $I\subset(1,R)$. For small $\delta$, \eqref{h-local-pattern} holds on this rectangle. Unfold the exact annular energy there. On every complete period square its indicator becomes exactly $\chi_h(z)$.

Replacing the smooth coefficients evaluated at $\delta a+\delta z$ by their values at $y$ changes the energy by $o(1)$, since the gradients are uniformly bounded in $L^2$. Omit squares cut by $\partial Q'$, whose energy is nonnegative. Weak lower semicontinuity and \eqref{h-radial-limit}--\eqref{h-angular-limit} give
\begin{align}\label{h-local-liminf}
\liminf_{\delta\to0}\int_{\Omega_\delta\cap X(I\times Q')}|\nabla v_\delta|^2\,dx
\geq{}&\int_{I\times Q'\times Y_h}J(y)
\Big[r^2|\partial_rw|^2\\
&\hspace{5mm}+ (\nabla_yw+\nabla_zW)^{\mathsf T}G(y)^{-1}
(\nabla_yw+\nabla_zW)\Big]\,dr\,dy\,dz.\notag
\end{align}
This follows first on rectangles a positive distance inside $Q'$, and then on $Q'$ by exhaustion; no bound on the discarded energy is needed.

At almost every $(r,y)$, the restriction of $W$ to $Y_h$ is an admissible periodic cell function. Minimizing its angular contribution in \eqref{h-local-liminf}, and integrating the radial contribution over $Y_h$, gives
\[
\int_{I\times Q'}J(y)
\Big[\theta r^2|\partial_rw|^2+B(G(y)^{-1},\nabla_yw)\Big]\,dr\,dy.
\]
Use \eqref{h-cell-A} and $dS=J(y)dy$, exhaust the interval and six open faces, and add the weak lower bound in $\B_1$. This proves
\[
\liminf_{\delta\to0}\q_\delta(v_\delta,v_\delta)\geq\q(v,v).
\]
For a finite lower limit, apply the extractions on a countable exhaustion to a subsequence realizing that limit.

Equation \eqref{h-radial-limit} proves, in particular, that the radial coefficient is exactly $\theta$. The original metric in \eqref{h-exact} has no radial--angular product, and its minimizing periodic adjustment only changes the angular derivative. Consequently the limiting energy has no mixed radial--angular term.

\paragraph{6. Approximating the energy from above.}
We first treat a smooth $v$ on $\overline\B_R$. Put $w=v\circ X$ on each face. For every $(r,y)$ choose the normalized minimizing cell function
\[
\zeta_*(r,y,z)=\zeta(G(y)^{-1},\nabla_yw(r,y)).
\]
It depends continuously on $(r,y)$ with values in $H^1(Y_h)$, by the estimate following \eqref{h-cell-A}, and is bounded there for this fixed $h$. Extend it into the full periodic square by the fixed-cell extension described above.

Cutting off in $(r,y)$ and then convolving, periodically in $z$, approximates these extensions in $L^2((1,R)\times Q;H^1_{\rm per}(Y))$ by smooth functions $W(r,y,z)$, periodic in $z$, whose support in $(r,y)$ is compactly contained in $(1,R)\times Q$. Dominated convergence justifies the cutoffs; derivatives in $(r,y)$ are not part of this norm. Since the quadratic integrals are continuous in $\nabla_zW$ for this norm,
\begin{equation}\label{h-approx-cell}
\sum_{\rm faces}\int_1^R\int_QJ(y)\int_{Y_h}
(\nabla_yw+\nabla_zW)^{\mathsf T}G^{-1}(\nabla_yw+\nabla_zW)\,dz\,dy\,dr
\end{equation}
can be made arbitrarily close to the limiting angular energy. Thus the exclusions near $r=1$ and the face boundaries have arbitrarily small energy cost before $\delta$ is sent to zero.

Keep such smooth functions $W$ fixed for the moment. Define a function on the full ball by
\begin{equation}\label{h-recovery}
V_\delta(x)=v(x)\quad(|x|\leq1),\qquad
V_\delta(X(r,y))=w(r,y)+\delta W(r,y,y/\delta)\quad(1<r<R).
\end{equation}
The second definition is made separately on each open face. Because $W$ vanishes near all face boundaries and near $r=1$, the definitions agree on their overlaps and traces. Thus $V_\delta\in H^1(\B_R)$, and $V_\delta\to v$ in $L^2(\B_R)$. Set $v_\delta=V_\delta|_{\Omega_\delta}$. Differentiation gives
\[
\partial_r(V_\delta\circ X)=\partial_rw+\delta\partial_rW(r,y,y/\delta),
\]
and
\[
\nabla_y(V_\delta\circ X)=\nabla_yw+
\nabla_zW(r,y,y/\delta)+\delta\nabla_yW(r,y,y/\delta).
\]
The terms multiplied by $\delta$ tend uniformly to zero for these fixed smooth $W$. Periodic averaging on complete period squares and \eqref{h-exact} now prove convergence of the radial energy to
$\theta\int r^2|\partial_rv|^2\,dr\,dS$,
and convergence of the angular energy to \eqref{h-approx-cell}. The additionally retained coordinate strips have area $O(\delta)$; since $W$ is zero near them, their energy is at most $C\delta$. The energy in the unit ball is exactly that of $v$.

Choose successive $W_n$ with angular energy error at most $1/n$, and then thresholds $d_n\downarrow0$ below which the fixed-$W_n$ convergence errors are at most $1/n$. Using $W_n$ when $d_{n+1}\leq\delta<d_n$ gives
\begin{equation}\label{h-recovery-conv}
\q_\delta(v_\delta,v_\delta)\longrightarrow\q(v,v),
\qquad
\|v_\delta-v\|_{L^2(\Omega_\delta)}\longrightarrow0.
\end{equation}
The extensions $E_\delta v_\delta$ are bounded in $H^1$. For any strong $L^2$ subsequential limit $a$, \eqref{h-mass-conv} applied to $E_\delta v_\delta-v$ gives
\[
\int_{\B_R}m|a-v|^2\,dx
=\lim\int_{\Omega_\delta}|v_\delta-v|^2\,dx=0.
\]
Since $m\geq\theta>0$, every such limit equals $v$, proving convergence of the full extended sequence. Smooth density in $H^1(\B_R)$, continuity of the limiting norm and energy, and a second diagonal choice give the recovery sequence for general $v\in H^1(\B_R)$.

\paragraph{7. Indexed eigenvalues.}
For fixed $R,h$, the limiting norm and energy are equivalent to the usual $L^2$ norm and gradient energy, by \eqref{coercive}. Compact embedding therefore gives a discrete spectrum and the min--max formula.

Take individual recovery sequences for the first $k+1$ limiting orthonormal eigenfunctions $e_0,\ldots,e_k$. Their mass matrices converge to the identity. Their energy matrices are bounded by Cauchy--Schwarz, and any subsequential limit $M$ satisfies, by the lower bound,
\[
a^{\mathsf T}Ma\geq\q\Big(\sum_i a_ie_i,\sum_i a_ie_i\Big)
\quad(a\in\R^{k+1}).
\]
The difference of these quadratic forms is positive semidefinite with zero diagonal, since each individual sequence recovers its energy, and hence is zero. Thus the energy matrices converge, without requiring linear recovery choices. The min--max formula on their spans gives
\[
\limsup_{\delta\to0}\lambda_k(\Omega_\delta)\leq\lambda_k.
\]

Conversely, along a subsequence realizing $\ell=\liminf\lambda_k(\Omega_\delta)$, extract simultaneous limits of the extensions of the first $k+1$ orthonormal actual eigenfunctions $u_{0,\delta},\ldots,u_{k,\delta}$. The upper bound bounds their energies, and \eqref{h-mass-conv} makes the limits $u_0,\ldots,u_k$ weighted orthonormal. Since
\[
\q_\delta\Big(\sum_i a_iu_{i,\delta},\sum_i a_iu_{i,\delta}\Big)
\leq\lambda_k(\Omega_\delta)\sum_i a_i^2,
\]
the lower energy bound and limiting min--max formula give $\lambda_k\leq\ell$. This proves indexed convergence with multiplicities.

\paragraph{8. The first eigenfunction and simplicity.}
Suppose $\lambda_1$ is simple. Then $\lambda_1<\lambda_2$, so indexed convergence implies
$\lambda_1(\Omega_\delta)<\lambda_2(\Omega_\delta)$
for all sufficiently small $\delta$. Thus the actual first positive eigenvalue is simple. Every strong $L^2$ limit of its normalized extended eigenfunctions has weighted norm one, weighted mean zero, and, by the lower bound, energy at most $\lambda_1$. It is therefore a minimizer for the first positive limiting eigenvalue. Simplicity makes it equal to $u_0$ or $-u_0$, where $u_0$ is any fixed normalized real limiting eigenfunction. Choose the sign so that
\[
\int_{\Omega_\delta}u_\delta u_0\,dx>0.
\]
This integral tends along each convergent subsequence to $\langle u,u_0\rangle_m=\pm1$, so it is nonzero for all sufficiently small $\delta$ and the chosen limit is always $u_0$. Compactness then proves convergence of the entire signed sequence. This completes every assertion of Lemma~\ref{hom}.\hfill$\square$

\subsection{Proof of Lemma~\ref{moser}}
\label{proofmoser}

Fix $R,\theta$. The extension lemma and the Sobolev inequality on the
fixed ball give, uniformly in $\delta$,
\begin{equation}\label{moserSobolev}
\|z\|_{L^6(\Omega_\delta)}^2
\leq S\Big(\|\nabla z\|_{L^2(\Omega_\delta)}^2
                  +\|z\|_{L^2(\Omega_\delta)}^2\Big),
\qquad z\in H^1(\Omega_\delta).
\end{equation}
All constants in the boundedness estimates below depend only on $S$
and, for the global estimate, the eigenvalue bound $\Lambda$.

Put $D_t=\Omega_\delta\cap B(x,t)$. The local weak formulation is
\begin{equation}\label{moserWeakLocal}
\int_{D_{2s}}\nabla v\cdot\nabla\varphi\,dx
=\int_{D_{2s}}F\varphi\,dx
\end{equation}
for every $\varphi\in H^1(D_{2s})$ whose support stays a positive
distance from $\partial B(x,2s)$. The test function may be nonzero
on $\partial\Omega_\delta$.

\paragraph{The local estimate.}
Take $s_0\leq1$, write $f=\|F\|_{L^\infty(D_{2s})}$, and set
$a=s^2f+\varepsilon$, $w=(v-c)_++a$, where $\varepsilon>0$.
Testing \eqref{moserWeakLocal} with an increasing approximation to
the indicator of $\{v>c\}$ times a nonnegative test function gives
$\int\nabla w\cdot\nabla\psi\leq f\int\psi$.
For $p\geq2$ and a Lipschitz cutoff $\eta$ supported strictly inside
$B(x,2s)$, testing this inequality with $\eta^2w^{p-1}$ and using
Young's inequality gives
\[
\int_{D_{2s}}|\nabla(\eta w^{p/2})|^2
\leq Cp^2\int_{D_{2s}}w^p
       \big(|\nabla\eta|^2+s^{-2}\eta^2\big),
\]
since $f/w\leq f/a\leq s^{-2}$. These tests do not assume
boundedness: use $\eta^2w\min(w,L)^{p-2}$ first and pass to the
limit as $L\to\infty$. The estimate is then valid whenever
$w\in L^p$.

Extend $\eta w^{p/2}$ by zero to the rest of $\Omega_\delta$ and
apply \eqref{moserSobolev}. Choosing $\eta=1$ on $B(x,t')$,
supported in $B(x,t)$, with $|\nabla\eta|\leq C/(t-t')$, yields
\[
\|w\|_{L^{3p}(D_{t'})}
\leq\big[CSp^2\big((t-t')^{-2}+s^{-2}\big)\big]^{1/p}
       \|w\|_{L^p(D_t)},
\qquad s\leq t'<t\leq2s.
\]
The endpoint $t=2s$ follows by approximation with cutoffs supported
strictly inside that ball. Iterate with $p_j=2\cdot3^j$ and
$t_j=s(1+2^{-j})$.
Since $\sum_jp_j^{-1}=3/4$ and $\sum_jj/p_j<\infty$, this gives
\[
\|w\|_{L^\infty(D_s)}
\leq Cs^{-3/2}\|w\|_{L^2(D_{2s})}
\leq C\big(s^{-3/2}\|v-c\|_{L^2(D_{2s})}+a\big),
\]
where the last inequality uses $|D_{2s}|\leq C s^3$.
Let $\varepsilon\downarrow0$ and repeat for $(c-v)_+$ to obtain
\eqref{local}. Only \eqref{moserSobolev} and Euclidean cutoffs were
used, so the constants remain uniform near every boundary piece
and its corners.

\paragraph{The global bound and boundary values.}
For a normalized Neumann eigenfunction with $0\leq\lambda\leq\Lambda$,
testing the equation with $|u|^{p-2}u$, justified by truncation as above,
and applying \eqref{moserSobolev} gives
\[
\|u\|_{L^{3p}(\Omega_\delta)}
\leq\big[CS p^2(1+\Lambda)\big]^{1/p}
       \|u\|_{L^p(\Omega_\delta)},\qquad p\geq2.
\]
Iteration starting at $p=2$ proves \eqref{globalbound}, with a
constant depending only on $R,\theta,\Lambda$. In particular,
$\lambda u$ is bounded. For each fixed $\delta$, the domain is
bounded and Lipschitz, so the Neumann boundary regularity theorem
\cite[Proposition~3.6]{Nittka} gives a continuous representative of
$u$ on $\overline\Omega_\delta$. No uniform boundary modulus of
continuity is needed. By continuity, the essential-supremum bound
\eqref{local} also holds at boundary points, which are limits of
interior points.\qed

\subsection{Proof of Lemma~\ref{trial}}\label{prooftrial}

In this proof only the constant weight $\theta$ used for the basic
starlike construction is considered.  Write
\[
 d\mu(r)=4\pi p(r)\,dr,
 \qquad
 p(r)=\begin{cases}r^2,&0<r<1,\\ \theta r^2,&1<r<R.\end{cases}
\]
The Rayleigh quotient on the radial subspace is
\[
 \frac{\int_0^R p(r)|v'(r)|^2\,dr}
      {\int_0^R p(r)|v(r)|^2\,dr},
 \qquad \int_0^R p(r)v(r)\,dr=0.
\]
The factors $4\pi$ cancel.  Define the continuous test function
\[
 g(r)=\begin{cases}
 1,&0\leq r\leq1,\\[2pt]
 \displaystyle\frac{1/r-1/R}{1-1/R},&1\leq r\leq R.
 \end{cases}
\]
Its radial extension belongs to $H^1(\B_R)$.  Its derivative has a jump only at $r=1$, where $g$ itself is continuous.  Direct
integration gives
\begin{align*}
 \int_0^R p(r)\,dr
   &=\frac{1+\theta(R^3-1)}{3},\\
 \int_0^R p(r)g(r)\,dr
   &=\frac13+\frac{\theta(R-1)(R+2)}6,\\
 \int_0^R p(r)|g'(r)|^2\,dr
   &=\frac{\theta}{1-1/R}.
\end{align*}
Let
\[
 b_R=\frac{\int_0^R p g}{\int_0^R p}
 =\frac{1+\theta(R-1)(R+2)/2}{1+\theta(R^3-1)}.
\]
For $\theta=R^{-5/2}/(4\pi)$, the numerator tends to $1$, whereas
the denominator is asymptotic to $R^{1/2}/(4\pi)$.  Thus $b_R\to0$.
For all sufficiently large $R$, we have $0<b_R\leq1/2$.

The function $g-b_R$ has the required mean zero.  On $[0,1]$ it
equals $1-b_R$, so
\[
 \int_0^R p(r)(g(r)-b_R)^2\,dr
 \geq\frac{(1-b_R)^2}{3}\geq\frac1{12}.
\]
Its energy is unchanged by subtracting the constant.  For $R\geq2$,
the variational characterization therefore yields
\[
 \lambda_{\rm rad}
 \leq\frac{\theta/(1-1/R)}{1/12}
 \leq24\theta
 =\frac6\pi R^{-5/2}.
\]
This proves \eqref{radbound} with a constant independent of $R$.
\hfill$\square$

\subsection{Proof of Lemma~\ref{sturm}}\label{proofsturm}

We first justify the variational eigenfunction and its boundary
conditions, including the condition at the singular endpoint $0$.
For the fixed $R$ and $\theta>0$ under consideration, the weighted
norm is equivalent to the ordinary $L^2(\B_R)$ norm, and the
weighted energy is equivalent to the ordinary Dirichlet energy.
The radial subspace of $H^1(\B_R)$ is closed.  The compact embedding
$H^1(\B_R)\to L^2(\B_R)$ and the Poincar\'e inequality on the ball
therefore give a minimizer $f$ for
\[
 \lambda_{\rm rad}
 =\inf\left\{
 \frac{\int_0^R p(v')^2}{\int_0^R pv^2}:
  v(|\cdot|)\in H^1(\B_R),\quad
  v\ne0,\quad \int_0^R pv=0
 \right\}.
\]
In particular $\lambda_{\rm rad}>0$.  A minimizing sequence may be
normalized by $\int_0^R pv^2=1$.  Its energies and its ordinary
$H^1$ norms are bounded; its weak $H^1$ limit is radial, and strong
$L^2$ convergence preserves both the normalization and the mean
condition.  Lower semicontinuity of the energy then gives the stated
minimizer.

Put $\lambda=\lambda_{\rm rad}$.  Differentiating the Rayleigh
quotient along mean-zero radial variations gives
\begin{equation}\label{radweak}
 \int_0^R p f'\varphi'\,dr
   =\lambda\int_0^R p f\varphi\,dr
\end{equation}
for every mean-zero radial test function.  Subtracting a constant
from an arbitrary radial test function does not change its derivative;
it also does not change the right side, since $\int_0^R pf=0$.
Consequently \eqref{radweak} holds for every radial test function.

On each of $(0,1)$ and $(1,R)$, distributional differentiation gives
\begin{equation}\label{radeq}
 -(pf')'=\lambda pf.
\end{equation}
On compact subintervals away from $0$, both $p$ and $p^{-1}$ are
bounded.  The equation implies that $pf'$ is absolutely continuous
there.  Since $p$ is smooth on each of the two intervals, the solution
is smooth there.  The function $f$ is continuous at $1$, because a
one-dimensional $H^1$ function on a neighborhood of $1$ is continuous.
Integrating \eqref{radweak} by parts on the two intervals, using test
functions supported near $1$, gives
\[
 p(1-)f'(1-)=p(1+)f'(1+).
\]
Using tests supported near $R$, with an arbitrary value at $R$, gives
$p(R)f'(R)=0$, hence $f'(R)=0$.

\textbf{The endpoint $0$.}
Let $k=\sqrt\lambda$.  In $(0,1)$ equation \eqref{radeq} is
$f''+2f'/r+\lambda f=0$.  Substituting $w(r)=rf(r)$ gives
$w''+\lambda w=0$.  Hence
\[
 f(r)=\frac{a\sin(kr)+b\cos(kr)}r,
 \qquad 0<r<1.
\]
If $b\ne0$, then $f'(r)=-b/r^2+O(1)$ as $r\downarrow0$, which
would give $\int_0^1 r^2|f'(r)|^2\,dr=\infty$.  The finite-energy
condition therefore forces $b=0$.  The remaining solution extends
smoothly across $0$ as a radial function and has
\begin{equation}\label{radregular}
 f(0)=ak,\qquad f'(0)=0,\qquad
 f(r)=f(0)\frac{\sin(kr)}{kr},\qquad
 \lim_{r\downarrow0}p(r)f'(r)=0.
\end{equation}
The coefficient $a$ is nonzero: otherwise $f$ vanishes in the unit
ball, and the value and flux at $1$ vanish.  Uniqueness for the
ordinary differential equation on $[1,R]$ would then make $f$
identically zero, a contradiction.  Thus $f(0)\ne0$, and its sign
can be chosen positive.

\textbf{Simplicity.}
Any radial eigenfunction with eigenvalue $\lambda>0$ satisfies the
same finite-energy condition at $0$.  By \eqref{radregular}, it is
uniquely determined on $[0,1]$ by its value at $0$.  Continuity and
flux continuity specify its value and derivative at $1+$, and
uniqueness for the differential equation then specifies it on
$[1,R]$.  Two radial eigenfunctions with eigenvalue $\lambda$ are
therefore proportional.  In particular the first positive radial
eigenvalue is simple.

\textbf{The number of zeros.}
There is at least one zero in $(0,R)$, since $f$ is nonzero and
$\int_0^R pf=0$.  A zero in either smooth interval is simple:
if both $f$ and $pf'$ vanished there, uniqueness for the equation
would imply $f\equiv0$.  The same argument works at $1$, using flux
continuity.  Also $f(R)\ne0$, since otherwise the Neumann condition
would give both value and derivative zero at $R$.  Formula
\eqref{radregular} excludes zeros close to $0$.  Zeros cannot
accumulate in the remainder of $[0,R]$, because an accumulation
point would have both value and derivative (or value and flux, at
$1$) zero.  Thus there are finitely many zeros, and the intervals
between consecutive zeros are precisely the intervals on which
$f$ has a fixed nonzero sign.

Suppose there were at least three such intervals.  Choose two of
them, $I_1$ and $I_2$, and let $f_j$ be $f$ on $I_j$ and zero
elsewhere.  Their radial extensions belong to $H^1(\B_R)$:
at each internal endpoint $f=0$, so the zero extension introduces
no jump.  Testing \eqref{radweak} with $f_j$ gives
\[
 \int_0^R p(f_j')^2\,dr
   =\lambda\int_0^R p f_j^2\,dr.
\]
Each number $m_j=\int_0^R p f_j\,dr$ is nonzero, because $f_j$
has a fixed nonzero sign on its interval.  Set
\[
 v=m_2f_1-m_1f_2.
\]
This is a nonzero radial function with $\int_0^R pv=0$.  Since
$f_1$ and $f_2$ have disjoint interiors of support, their mass and
energy cross terms vanish.  Therefore
\[
 \int_0^R p(v')^2\,dr
   =\lambda\int_0^R pv^2\,dr.
\]
Thus $v$ is another minimizer for the first positive radial
Rayleigh quotient, and the variational argument above shows that
it is a radial eigenfunction with eigenvalue $\lambda$.
Simplicity makes $v$ proportional to $f$.  But $v$ vanishes on
the third sign interval, where $f$ does not vanish.  This is
impossible.  There are therefore exactly two sign intervals,
separated by exactly one zero $z\in(0,R)$.

\textbf{Strict decrease.}
The chosen sign gives $f>0$ on $[0,z)$ and $f<0$ on $(z,R]$.
Define $J(r)=\int_0^r p(s)f(s)\,ds$.  Then $J(r)>0$ for
$0<r\leq z$.  If $z<r<R$, the mean-zero condition gives
\[
 J(r)=-\int_r^R p(s)f(s)\,ds>0.
\]
Integrating \eqref{radeq}, using \eqref{radregular} and the flux
condition at $1$, gives
\[
 p(r)f'(r)=-\lambda J(r)<0\qquad(0<r<R),
\]
with the derivative interpreted one-sidedly at $1$.  The function
is therefore strictly decreasing on $[0,R]$.
This proves every assertion of Lemma~\ref{sturm}, and also proves
the monotonicity used in the main argument.
\hfill$\square$

\end{document}